\documentclass[a4paper, 12pt]{article}

\usepackage{epsfig, color, epstopdf}
\usepackage{amssymb,amsmath,amsthm,latexsym,lineno}
\usepackage{hyperref}
\usepackage{graphicx,graphics,amsfonts,indentfirst}
\usepackage{mathrsfs,eufrak, bbm, dsfont,yfonts}
\usepackage{url,tikz,yhmath,appendix,booktabs}
\usepackage{multirow}
\usepackage{kantlipsum}
\allowdisplaybreaks
\usepackage{setspace}

\newtheorem{theo}{Theorem}
\newtheorem{lem}[theo]{Lemma}
\newtheorem{pro}[theo]{Proposition}
\newtheorem{prob}{Problem}
\newtheorem{exa}[theo]{Example}
\newtheorem{con}{Conjecture}
\newtheorem{cor}[theo]{Corollary}

\usetikzlibrary{backgrounds}
\usetikzlibrary{topaths,calc}
\usetikzlibrary{arrows}
\usetikzlibrary{shapes,shapes.geometric,shapes.misc}
\pgfdeclarelayer{edgelayer}
\pgfdeclarelayer{nodelayer}
\pgfsetlayers{background,edgelayer,nodelayer,main}
\tikzstyle{none}=[inner sep=0mm]
\tikzstyle{every loop}=[]

\tikzstyle{dotted}=[dash pattern=on \pgflinewidth off 2pt]
\tikzstyle{dashed}=[dash pattern=on 3pt off 3pt]

\newcommand \tikzp[2]
{
	\begin{center}
		\begin{tikzpicture}[scale=#1]
			#2
		\end{tikzpicture}
	\end{center}
}

\tikzstyle{new style 0}=[fill=black, draw=black, shape=circle]
\tikzstyle{red style 1}=[fill=red, draw=black, shape=circle]
\tikzstyle{blue style 2}=[fill=blue, draw=black, shape=circle]
\tikzstyle{white style 4}=[fill=white, draw=black, shape=circle]
\tikzstyle{bklack style 5}=[fill=black, draw=black, shape=rectangle]
\tikzstyle{red style 3}=[fill=red, draw=black, shape=rectangle]
\tikzstyle{yellow style 7}=[fill=yellow, draw=black, shape=rectangle]
\tikzstyle{new style 8}=[fill={rgb,255: red,0; green,132; blue,0}, draw={rgb,255: red,0; green,131; blue,0}, shape=circle]

\tikzstyle{new edge style 0}=[-]
\tikzstyle{new edge style 1}=[-, draw=red]
\tikzstyle{new edge style 2}=[-, draw=blue]
\tikzstyle{new edge style 3}=[-, draw={rgb,255: red,0; green,156; blue,0}]

\tikzstyle{cblue}=[circle, draw, thin,fill=blue!20, scale=0.5]

\newcommand \equ[2]
{
	\begin{equation}\label{#1}
		#2
	\end{equation}
}

\newcommand \aln[2]
{
	\begin{align}\label{#1}
		#2
	\end{align}
}

\newcommand \eqn[2]
{
	\begin{eqnarray}\label{#1}
		#2
	\end{eqnarray}
}

\newcommand \lemm[2]
{\begin{lem}
		\label{#1} #2
	\end{lem}
}

\newcommand \prop[2]
{\begin{pro}
		\label{#1} #2
	\end{pro}
}

\newcommand \corr[2]
{\begin{cor}
		\label{#1} #2
	\end{cor}
}

\newcounter{countcase}

\newcounter{countclaim}

\def \proof {\noindent {{\it Proof}}.\setcounter{countcase}{0} \setcounter{clm}{0}}

\newcommand{\proofend}{{\hfill$\Box$}\setcounter{countclaim}{0}\setcounter{countcase}{0}}

\def \N {{\mathbb N}}

\def \hyh {{\cal H}}
\newcommand \spann[1]{\langle #1\rangle }

\def \nbroken {{\mathscr {N}}\hspace{-0.15 cm}{\mathscr {B}}}

\begin{document}
	\baselineskip 0.6 cm
	
	\title{
		Compare list-color functions of uniform
		hypergraphs	with their chromatic polynomials (II)
}
	
\author{
Meiqiao Zhang\thanks{Corresponding Author. Email: nie21.zm@e.ntu.edu.sg and 
								meiqiaozhang95@163.com.},
		Fengming Dong\thanks{
			Email: fengming.dong@nie.edu.sg and donggraph@163.com.}
		\\
		National Institute of Education,
		Nanyang Technological University, Singapore
	}
	\date{}
	
	\maketitle{}

	\begin{abstract}
	For any $r$-uniform hypergraph $\mathcal{H}$ with $m$ ($\geq 2$) edges,
let $P(\mathcal{H},k)$ and $P_l(\mathcal{H},k)$ 
 be the chromatic polynomial and 
 the list-color function of $\mathcal{H}$
 respectively,
 and let 
 $\rho(\mathcal{H})$ denote the minimum value 
 of $|e\setminus e'|$ among 
 all pairs of distinct edges $e,e'$ in $\mathcal{H}$. 
We will show that if $r\ge3$, $\rho(\mathcal{H})\ge 2$ and $m\ge \frac{\rho(\mathcal{H})^3}2+1$, then
$P_l(\mathcal{H},k)=P(\mathcal{H},k)$ holds 
for all integers $k\geq \frac{2.4(m-1)}{\rho(\mathcal{H})\log(m-1)}$.
	\end{abstract}



\section{ Introduction
\label{secintro}}
In this article, we always assume that $\hyh$ is a hypergraph with vertex set $V(\hyh)$ and edge set $E(\hyh)$, where $|e|\ge 2$ for every edge $e$ in $\hyh$ and $e_1\not\subseteq e_2$ for any pair of edges $e_1,e_2$ in $\hyh$. A hypergraph $\hyh$ is called \textit{$r$-uniform} if $|e|=r$ for every $e\in E(\hyh)$, and $\hyh$ is called \textit{linear} if $|e_1\cap e_2|\le 1$ 
	for any two edges $e_1,e_2$ in $\hyh$. Then obviously, any simple graph is $2$-uniform and linear. 	
The number of connected components of $\hyh$ is denoted by $c(\hyh)$.
	For any edge set $F\subseteq E(\hyh)$, 
	let $\hyh\spann{F}$ be the spanning subhypergraph of $\hyh$ with edge set $F$, let $c(F)=c(\hyh\spann{F})$, 
	and let $V(F)$ be the set of vertices 
	$v$ in $\hyh$ with $v\in e$ for some $e\in F$.
Let $\N$ denote the set of positive integers and $[k]=\{1,2,\dots,k\}$ for any $k\in\N$.

For any $k\in\N$, a \textit{proper $k$-coloring} of $\hyh$ is a mapping $\theta:V(\hyh)\rightarrow [k]$ such that $|\{\theta(v):   v\in e\}|\ge 2$ for each edge $e$ in $\hyh$. Then, as a generalization of proper coloring, \textit{list-coloring} was introduced by Vizing~\cite{vizing}, and Erd\H{o}s, Rubin and Taylor~\cite{erdos} independently.
For any $k\in\N$, a \textit{$k$-assignment} $L$ of $\hyh$ is a mapping from $V(\hyh)$ to the power set of $\N$ such that $|L(v)|=k$ for each vertex $v$ in $\hyh$, and an \textit{$L$-coloring} of $\hyh$ is a mapping $\theta:V(\hyh)\rightarrow \N$ such that $\theta(v)\in L(v)$ for each $v\in V(\hyh)$ and 
$|\{\theta(v): v\in e\}|\ge 2$ for each 
$e\in E(\hyh)$. 

From a counting point of view, the \textit{chromatic polynomial} $P(\hyh, k)$ of $\hyh$ is defined to be the number of proper $k$-colorings of $\hyh$ for any $k\in\N$. Similarly, for any $k$-assignment $L$ of $\hyh$, let $P(\hyh, L)$ be the number of $L$-colorings of $\hyh$, and the \textit{list-color function} $P_l(\hyh,k)$ of $\hyh$ counts the minimum value of $P(\hyh, L)$'s among all $k$-assignments $L$ of $\hyh$ for any $k\in\N$. Then, $P_l(\hyh,k)\le P(\hyh,k)$ follows directly, and Kostochka and Sidorenko~\cite{Kosto} raised the question that under what conditions, the equality holds.
For the history of 
the study on this problem on graphs
and uniform hypergraphs,
we refer our readers to \cite{Dong22b, Dong22, Thomassen, wang20, wang17}. 
In particular,  this article is a continuation of 
\cite{Dong22b}, focusing on uniform hypergraphs, in which 
$|e\setminus e'|\ge 2$ for each pair of 
distinct edges $e$ and $e'$.
More information on these two color functions can be found in~\cite{birk, dong1, dong0, Jack15, Kaul22, rea1, rea2, Tomescu1998, Tomescu2009,ruixue1,ruixue2}.  

In what follows, assume that 
$\hyh$ is an $r$-uniform hypergraph with $m$ edges,
where $r\ge 3$. 
Let $E_{r-1}(e)$ be the set of edges $e'$ 
in $\hyh$ with $|e\cap e'|=r-1$,
and let 
$\gamma(\hyh)$ be the maximum size 
of $E_{r-1}(e)$ among all edges 
$e$ in $\hyh$. 
In \cite{Dong22b}, the authors have shown that 
$P_l(\hyh,k)=P(\hyh,k)$ holds 
whenever $k\ge \min\{m-1,0.6(m-1)+0.5\gamma(\hyh)\}$,
which improved the result due to Wang, Qian and Yan~\cite{wang20}
that $P_l(\hyh,k)=P(\hyh,k)$ holds 
whenever $k\ge 1.1346(m-1)$. 
In this article, we shall further improve this conclusion for the case when $\gamma(\hyh)=0$.

Since $P_l(\hyh,k)=P(\hyh,k)$ trivially holds whenever $m\le 1$ and $k\in\N$, we assume that $m\ge 2$ in the following.
Let $\rho(\hyh):=\min\{|e\setminus e'|:
e\ne e', e,e'\in E(\hyh)\}$.
Clearly, $1\le \rho(\hyh)\le r$, 
$\hyh$ is linear if and only if $\rho(\hyh)\ge r-1$,
and $\gamma(\hyh)=0$ if and only if $\rho(\hyh)\ge 2$. Then our main result is stated as follows.

 \begin{theo}\label{th4-4}
 Let $\hyh$ be an $r$-uniform hypergraph with 
 $m$ edges, where $r\ge 3$. 
If $\rho(\hyh)\ge 2$ and $m\ge \rho(\hyh)^3/2+1$,
then $P_l(\hyh,k)=P(\hyh,k)$ holds 
whenever $k\ge 
\frac{2.4(m-1)}{\rho(\hyh)\log(m-1)}$. 
 \end{theo}

Note that Theorem~\ref{th4-4}  
presents a result much better than 
that in \cite{Dong22b} under the condition $\rho(\hyh)\ge 2$ and $m\ge \rho(\hyh)^3/2+1$.
In the following, we shall give more detailed comparison results for linear and uniform hypergraphs.

For any $k$-assignment $L$ of $\hyh$, let
$\alpha(\hyh, L):=\sum_{e\in E(\hyh)}\alpha(e,L)$, 
where 
$\alpha(e,L):=k-\left|\bigcap\limits_{v\in e}L(v)\right|$ for each $e\in E(\hyh)$.
Obviously, $0\le \alpha(e,L)\le k$, and 
$\alpha(e,L)=0$ if and only if 
$\left|\bigcup_{v\in e}L(v)\right|=k$
(i.e., $L(v)$ is a fixed $k$-element subset of $\N$ for all $v\in e$).
It implies that $\alpha(\hyh,L)=0$ 
if and only if $\left|\bigcup_{v\in e}L(v)\right|=k$
for each edge $e\in E(\hyh)$.
Further, it can be verified easily that 
$P(\hyh,L)=P(\hyh,k)$ holds when $\alpha(\hyh,L)=0$. Therefore, the obstacle here is the case when $\alpha(\hyh,L)>0$.

\begin{theo}\label{th4-1}
Let $\hyh$ be a linear and $3$-uniform hypergraph with $n$ vertices and $m$ ($\ge 3$) edges.
If $k\ge 
  \frac{1.185(m-1)}{\log(m-1)}$
  and $L$ is a $k$-assignment of $\hyh$ with $\alpha(\hyh,L)>0$, then
\aln{eq4-1}
{
P(\hyh,L)-P(\hyh,k)
>
\frac{0.002	\log(m-1)}
{(m-1)^{0.156}}
	\cdot 
k^{n-3}\alpha(\hyh,L).
}
Hence $P_l(\hyh,k)=P(\hyh,k)$ holds whenever $k\ge 
  \frac{1.185(m-1)}{\log(m-1)}$.
\end{theo}

 \begin{theo}\label{th4-3}
 Let $\hyh$ be a linear and $r$-uniform hypergraph with $n$ vertices and $m$ ($\ge 3$) edges, where $r\ge 4$. 
 If $k\ge 
    \frac{0.831(m-1)}{\log(m-1)}$ and 
$L$ is a $k$-assignment of $\hyh$ with $\alpha(\hyh,L)>0$, then
\aln{eq4-3}
{
P(\hyh,L)-P(\hyh,k)
>(m-1)^{-1.796}
 (1+1.796\log(m-1))k^{n-r}
 \alpha(\hyh,L).
}
Hence $P_l(\hyh,k)=P(\hyh,k)$ holds whenever $k\ge 
    \frac{0.831(m-1)}{\log(m-1)}$.
 \end{theo}

We will introduce some basic results in Section~\ref{sec2}, and then give the proofs of Theorems~\ref{th4-4},~\ref{th4-1} and~\ref{th4-3} in Section~\ref{sec5}.

\section{Preliminary results
\label{sec2}}
In this section, we always assume that $\hyh$ is a hypergraph with vertex set $V$ and edge set $E$, where $|V|=n$ and $|E|=m$, $L$ is any $k$-assignment of $\hyh$, and $\eta$ is a fixed bijection from $E$ to $[m]$. 
We shall first introduce some known generalizations~\cite{Trinks14,wang20} of Whitney's broken-cycle theorem in terms of hypergraphs
and then list some related results for later use.

\begin{figure}[!ht]
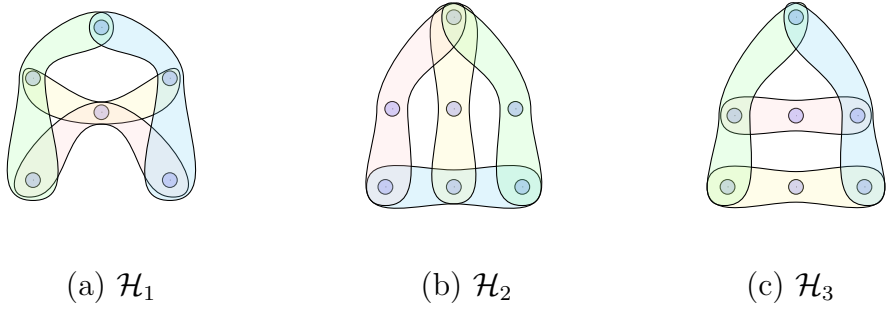

	\tikzp{0.45}
	{
		\foreach \place/\y in {{(-2.3,0.2)/1}, {(-0.3,1.2)/2},{(-4.3,1.2)/3},
			{(-0.3,-1.8)/4}, {(-4.3,-1.8)/5}, {(-2.3,2.7)/6}}
		\node[cblue] (b\y) at \place {};
		
		\filldraw[black] (b1) circle (0pt)node[anchor=south] {};
		\filldraw[black] (b2) circle (0pt)node[anchor=south] {};
		\filldraw[black] (b3) circle (0pt)node[anchor=south] {};
		\filldraw[black] (b4) circle (0pt)node[anchor=south] {};
		\filldraw[black] (b5) circle (0pt)node[anchor=north] {};
		\filldraw[black] (b6) circle (0pt)node[anchor=north] {};

    \begin{scope}[fill opacity=0.2]
    \filldraw[fill=yellow!50] ($(b1)+(0,0.3)$) 
        to[out=0,in=90] ($(b2) + (0.3,0)$) 
        to[out=270,in=0] ($(b1) + (0,-0.35)$)
        to[out=180,in=270] ($(b3) + (-0.3,0)$)
        to[out=90,in=180] ($(b1)+(0,0.3)$);  
    \filldraw[fill=pink!80] ($(b1)+(0,0.3)$) 
        to[out=0,in=0] ($(b4) + (0.1,-0.5)$) 
        to[out=180,in=0] ($(b1) + (0,-0.35)$)
        to[out=180,in=0] ($(b5) + (-0.1,-0.5)$)
        to[out=180,in=180] ($(b1)+(0,0.3)$);  
    \filldraw[fill=green!40] ($(b3)+(-0.4,0)$) 
        to[out=270,in=180] ($(b5) + (0,-0.6)$) 
        to[out=0,in=270] ($(b3) + (0.35,0)$)
        to[out=90,in=270] ($(b6) + (0.35,0)$)
        to[out=90,in=90] ($(b3)+(-0.4,0)$);          
    \filldraw[fill=cyan!50] ($(b2)+(-0.35,0)$) 
        to[out=270,in=180] ($(b4) + (0,-0.6)$) 
        to[out=0,in=270] ($(b2) + (0.4,0)$)
        to[out=90,in=90] ($(b6) + (-0.35,0)$)
        to[out=270,in=90] ($(b2)+(-0.35,0)$);             
    \end{scope}

		\foreach \place/\y in {{(8,3)/1}, {(6.2,0.3)/2},{(8,0.3)/3},
			{(9.8,0.3)/4}, {(6,-2)/5}, {(8,-2)/6},{(10,-2)/7}}
		\node[cblue] (b\y) at \place {};
		
		\filldraw[black] (b1) circle (0pt)node[anchor=south] {};
		\filldraw[black] (b2) circle (0pt)node[anchor=south] {};
		\filldraw[black] (b3) circle (0pt)node[anchor=south] {};
		\filldraw[black] (b4) circle (0pt)node[anchor=south] {};
		\filldraw[black] (b5) circle (0pt)node[anchor=north] {};
		\filldraw[black] (b6) circle (0pt)node[anchor=north] {};
		\filldraw[black] (b7) circle (0pt)node[anchor=north] {};
				
    \begin{scope}[fill opacity=0.2]
\filldraw[fill=cyan!50] ($(b5)+(-0.55,0)$) 
           to[out=90,in=180] ($(b6) + (0,0.5)$) 
           to[out=0,in=90] ($(b7) + (0.55,0)$)
           to[out=270,in=0] ($(b6) + (0,-0.5)$)
           to[out=180,in=270] ($(b5)+(-0.55,0)$);
\filldraw[fill=yellow!50] ($(b1)+(0,0.4)$) 
           to[out=0,in=90] ($(b3) + (0.45,0)$) 
           to[out=270,in=0] ($(b6) + (0,-0.5)$)
           to[out=180,in=270] ($(b3) + (-0.45,0)$)
           to[out=90,in=180] ($(b1)+(0,0.4)$);           
\filldraw[fill=pink!80] ($(b2)+(-0.45,0)$) 
           to[out=270,in=180] ($(b5) + (0,-0.55)$) 
           to[out=0,in=270] ($(b2) + (0.45,0)$)
           to[out=90,in=270] ($(b1) + (0.35,0)$)
           to[out=90,in=90] ($(b2)+(-0.45,0)$);  
 \filldraw[fill=green!40] ($(b4)+(-0.45,0)$) 
            to[out=270,in=180] ($(b7) + (0,-0.55)$) 
            to[out=0,in=270] ($(b4) + (0.45,0)$)
            to[out=90,in=90] ($(b1) + (-0.35,0)$)
            to[out=270,in=90] ($(b4)+(-0.45,0)$);                          
    \end{scope}    
    
    		\foreach \place/\y in {{(18,3)/1},{(16.2,0.1)/2},{(18,0.1)/3},
    		{(19.8,0.1)/4}, {(16,-2)/5},{(18,-2)/6},{(20,-2)/7}}
    		\node[cblue] (b\y) at \place {};
    		
    		\filldraw[black] (b1) circle (0pt)node[anchor=south] {};
    		\filldraw[black] (b2) circle (0pt)node[anchor=south] {};
    		\filldraw[black] (b3) circle (0pt)node[anchor=south] {};
    		\filldraw[black] (b4) circle (0pt)node[anchor=south] {};
    		\filldraw[black] (b5) circle (0pt)node[anchor=north] {};
    		\filldraw[black] (b6) circle (0pt)node[anchor=north] {};
    		\filldraw[black] (b7) circle (0pt)node[anchor=north] {};
    				
        \begin{scope}[fill opacity=0.2]
        \filldraw[fill=pink!80] ($(b2)+(-0.4,0)$) 
           to[out=90,in=180] ($(b3) + (0,0.45)$) 
           to[out=0,in=90] ($(b4) + (0.4,0)$)
           to[out=270,in=0] ($(b3) + (0,-0.45)$)
           to[out=180,in=270] ($(b2)+(-0.4,0)$);        \filldraw[fill=yellow!50] ($(b5)+(-0.6,0)$) 
           to[out=90,in=180] ($(b6) + (0,0.45)$) 
           to[out=0,in=90] ($(b7) + (0.6,0)$)
           to[out=270,in=0] ($(b6) + (0,-0.45)$)
           to[out=180,in=270] ($(b5)+(-0.6,0)$);
\filldraw[fill=green!40] ($(b2)+(-0.45,0)$) 
           to[out=270,in=180] ($(b5) + (-0.09,-0.55)$) 
           to[out=0,in=270] ($(b2) + (0.45,0)$)
           to[out=90,in=270] ($(b1) + (0.35,0)$)
           to[out=90,in=90] ($(b2)+(-0.45,0)$);  
 \filldraw[fill=cyan!50] ($(b4)+(-0.45,0)$) 
            to[out=270,in=180] ($(b7) + (0.09,-0.55)$) 
            to[out=0,in=270] ($(b4) + (0.45,0)$)
            to[out=90,in=90] ($(b1) + (-0.35,0)$)
            to[out=270,in=90] ($(b4)+(-0.45,0)$);                            
        \end{scope}
    	}
    {}\hfill \hspace{0.2cm} (a) $\hyh_1$ \hspace{3.2 cm} (b) $\hyh_2$ \hspace{2.8cm}  (c) $\hyh_3$  \hfill {}    	
	\caption{Three linear and $3$-uniform hypergraphs of size four}
	\label{fig1}
\end{figure}

In 2014, Trinks~\cite{Trinks14} defined a \textit{$\delta$-cycle} of a hypergraph to be a minimal edge set $F$ such that $e\subseteq V(F\setminus \{e\})$ for every $e\in F$. By definition, a $\delta$-cycle contains at least three edges, and in Figure~\ref{fig1}, $E(\hyh_1)$ is a $\delta$-cycle while $E(\hyh_2)$ and $E(\hyh_3)$ are not.
An edge set $B\subseteq E$ is called a \textit{broken-$\delta$-cycle} if $B$ is obtained from a $\delta$-cycle $C$ by deleting the edge $e$ in $C$ with $\eta(e)\le \eta(e')$ for all $e'\in C$. Then, a broken-$\delta$-cycle includes at least two edges.

Let $\nbroken(\hyh)$ be the set of subsets of $E$ containing no broken-$\delta$-cycles. Clearly, $\{e\}\in \nbroken(\hyh)$ for every edge $e\in E$.
Further, Trinks~\cite{Trinks14}  generalized Whitney's broken-cycle theorem on  
$P(G,k)$ for a simple graph $G$
(see \cite{Whitney1932})  
to  $P(\hyh,k)$ for a hypergraph $\hyh$:
\equ{thhyBC}
{
P(\hyh,k)=\sum_{A\in \nbroken(\hyh)}(-1)^{|A|}k^{c(A)}.
}

For any $A\subseteq E$, let 
$A_1,\dots,A_{c(A)}$ be the components of $\hyh\langle A\rangle$,
$$
\beta(A_i,L)=\left |\bigcap\limits_{v\in V(A_i)}L(v)\right | \quad\text{and}\quad
 \beta(A,L)=\prod\limits_{i=1}^{c(A)}\beta(A_i,L).
$$
A result on $P(\hyh, L)$ analogous to  (\ref{thhyBC})
was established by
Wang, Qian and Yan~\cite{wang20}: 
\equ{thhyBCL}
{
P(\hyh,L)=\sum_{A\in \nbroken(\hyh)}(-1)^{|A|}\beta(A,L).
}

By (\ref{thhyBC}) and (\ref{thhyBCL}), we have
\eqn{eq3-1}
{
P(\hyh,L)-P(\hyh,k)=\sum_{A\in \nbroken(\hyh)}(-1)^{|A|}\left(\beta(A,L)-k^{c(A)}\right).
}
Clearly, for any $A\subseteq E$, $\beta(A,L)-k^{c(A)}\le 0$ follows from the definition. On the other hand, a lower bound of $\beta(A,L)-k^{c(A)}$ was determined by Wang, Qian and Yan~\cite{wang20}.

\begin{lem}[\cite{wang20}]\label{th3-1}
Let $L$ be a $k$-assignment of $\hyh$.
Then for any $A\subseteq E$,
$$
\beta(A,L)-k^{c(A)}\ge -k^{c(A)-1}\sum_{e\in A}\alpha(e,L).
$$
\end{lem}

For any $e\in E$, let $\nbroken(\hyh,e)$ be the set of $A\in \nbroken(\hyh)$ 
with $e\in A$. Then based on Lemma~\ref{th3-1}, a rough comparison between $P(\hyh,L)$ and $P(\hyh,k)$ can be done.

\prop{pp1-2}
{
Let $\hyh$ be an $r$-uniform hypergraph with $n$ vertices and $L$ be a $k$-assignment of $\hyh$.
Then
\equ{eq3-8}
{
P(\hyh, L)-P(\hyh,k)
\ge
\sum_{e\in E}\alpha(e,L)\left(
k^{n-r}-\left(\sum_{A\in \nbroken(\hyh,e)\atop |A|~\text{even}}k^{c(A)-1}\right)
\right).
}
}

\proof
It is easy to see that when $E=\emptyset$, $P(\hyh,L)=P(\hyh,k)$ holds, and hence the result is proven. In the following, assume that $|E|\ge 1$.
Then by Lemma~\ref{th3-1}, (\ref{eq3-1}) can be transferred to 
\eqn{eq3-4}
{
P(\hyh,L)-P(\hyh,k)&=&\sum_{A\in \nbroken(\hyh)}(-1)^{|A|}\left(\beta(A,L)-k^{c(A)}\right)
\nonumber\\
&=&
\sum_{A\in \nbroken(\hyh)\atop |A|~\text{odd}}\left(k^{c(A)}-\beta(A,L)\right)
+\sum_{A\in \nbroken(\hyh)\atop |A|~\text{even}}\left(\beta(A,L)-k^{c(A)}\right)
\nonumber\\
&\ge&
\sum_{A\in \nbroken(\hyh)\atop |A|=1}
\left(k^{c(A)}-\beta(A,L)\right)+
\sum_{A\in \nbroken(\hyh)\atop |A|~\text{even}}
\left(-k^{c(A)-1}\sum_{e\in A}\alpha(e,L)\right)
\nonumber\\
&=&
\sum_{|A|=1}
\left(k^{c(A)}-\beta(A,L)\right)+
\sum_{A\in \nbroken(\hyh)\atop |A|~\text{even}}
\left(-k^{c(A)-1}\sum_{e\in A}\alpha(e,L)\right)
\nonumber\\
&=&
\sum_{e\in E}\alpha(e,L)k^{n-r}+
\sum_{e\in E}\left(\alpha(e,L)\sum_{A\in \nbroken(\hyh,e)\atop |A|~\text{even}}-k^{c(A)-1}\right)
\nonumber\\
&=&
\sum_{e\in E}\alpha(e,L)\left(
k^{n-r}-\left(\sum_{A\in \nbroken(\hyh,e)\atop |A|~\text{even}}k^{c(A)-1}\right)
\right).
}
\proofend

To apply Proposition~\ref{pp1-2} more efficiently, we shall analyze the number
$c(A)$ for any $A\in\nbroken(\hyh)$ next. The idea of the two following lemmas are based on Theorem 5 of~\cite{wang20}.

\lemm{le2-1}
{
Let $\hyh$ be an $r$-uniform hypergraph with $n$ vertices, where $r\ge 3$. For any $A\in\nbroken(\hyh)$, 
$c(A)= n-r+1$ if $|A|=1$,
and $c(A)\le n-r-\rho(\hyh)-|A|+3$ if $|A|\ge 2$.
}
\proof 
It is easy to verify that the result holds
when $1\le |A|\le 2$.
Now, assume that $i=|A|\ge 3$.
For any $F\in\nbroken(\hyh)$, 
$F$ does not include any $\delta$-cycle, 
implying that 
$e\not\subseteq V(F\setminus \{e\})$
for some edge $e\in F$.
Thus, 
the edges in $A$ can be labeled as
$e_1,e_2,\dots,e_i$ such that
\equ{}
{
c(\{e_1\})>c(\{e_1,e_2\})>c(\{e_1,e_2,e_3\})>\cdots>c(\{e_1,e_2,\dots,e_i\})=c(A).
}
As $c(\{e_1,e_2\})\le n-r-\rho(\hyh)+1$, we have
\equ{}
{
c(A)
\le (n-r-\rho(\hyh)+1)-(i-2)=n-i-r-\rho(\hyh)+3.
}
\proofend

In particular, for linear and uniform hypergraphs, we have a refining result.

\lemm{le2-2}
{Let $\hyh$ be a linear and $r$-uniform hypergraph with $n$ vertices, where $r\ge 3$. For any $A\in\nbroken(\hyh)$, 
$c(A)= n-(r-1)|A|$ if $1\le |A|\le 2$,
$n-3r+3\le c(A)\le n-3r+4$
if $|A|=3$, 
and $c(A)\le n-3r-|A|+7$ if $|A|\ge 4$.
}
\proof 
It is easy to verify that the result holds
when $1\le |A|\le 3$.
Now, assume that $i=|A|\ge 4$.
Then similarly, 
the edges in $A$ can be labeled as
$e_1,e_2,\dots,e_i$ such that
\equ{}
{
c(\{e_1\})>c(\{e_1,e_2\})>c(\{e_1,e_2,e_3\})>\cdots>c(\{e_1,e_2,\dots,e_i\})=c(A).
}
As $c(\{e_1,e_2,e_3\})\le n-3r+4$, we have
\equ{}
{
c(A)
\le (n-3r+4)-(i-3)=n-i-3r+7.
}
\proofend

For any $i\in[m]$, let $\nbroken_i(\hyh,e)$ be the set of $A\in \nbroken(\hyh,e)$ with $|A|=i$. Obviously, $|\nbroken_i(\hyh,e)|\le \binom{m-1}{i-1}$.
Then, the next two corollaries immediately follow from Proposition~\ref{pp1-2} and Lemmas~\ref{le2-1},~\ref{le2-2}.

\corr{co1-3}
{
Let $\hyh$ be an $r$-uniform hypergraph with $n$ vertices and $m$ edges, where $r\ge 3$. If $L$ is a $k$-assignment of $\hyh$ with $\alpha(\hyh,L)>0$, then
\equ{eq3-12}
{
\frac{P(\hyh, L)-P(\hyh,k)}
{k^{n-r} \alpha(\hyh, L)}
\ge
1-\sum_{i\ge 1} \binom{m-1}{2i-1}k^{-2i-\rho(\hyh)+2}
>
1-\sum_{i\ge 1} \frac{(m-1)^{2i-1}}{(2i-1)!}k^{-2i-\rho(\hyh)+2}.
}
}

\corr{co1-2}
{
Let $\hyh$ be a linear and $r$-uniform hypergraph with $n$ vertices and $m$ edges, where $r\ge 3$. If $L$ is a $k$-assignment of $\hyh$ with $\alpha(\hyh,L)>0$, then
\eqn{eq6-1}
{
	\frac{P(\hyh, L)-P(\hyh,k)}
	{k^{n-r}\alpha(\hyh, L)}
&\ge&
		1
		-(m-1)k^{-r+1}
		-\sum_{i\ge 2}
	\binom{m-1}{2i-1}k^{-2i-2r+6}
\nonumber\\
&>&
	1
	-(m-1)(k^{-r+1}-k^{-2r+4})
	-\sum_{i\ge 1}
\frac{(m-1)^{2i-1}}{(2i-1)!}k^{-2i-2r+6}.
}
}

\section{Proofs of Theorems
\label{sec5}}

In the section, we shall prove our main results. 


 \noindent\textbf{Proof of Theorem~\ref{th4-4}.}
We need only to show that under the given conditions, $P(\hyh,L)>P(\hyh,k)$ holds for any $k$-assignment $L$ of $\hyh$ with $\alpha(\hyh,L)>0$.
 
Let $n:=|V(\hyh)|$, $t:=\rho(\hyh)$, $M:=m-1$ and $L$ be a $k$-assignment of $\hyh$ with $\alpha(\hyh,L)>0$. Then by (\ref{eq3-12}),
 \eqn{eq3-70}
 {
 \frac{P(\hyh, L)-P(\hyh,k)}
 {k^{n-r} \alpha(\hyh, L)}
 &>&
 1-\sum_{i\ge 1} \frac{M^{2i-1}}{(2i-1)!}k^{-2i-t+2}
 \nonumber\\
 &=&
1- k^{1-t}\left(\frac{\exp(M/k)-\exp(-M/k)}{2}\right)
\nonumber\\
&>&
1- \frac{k^{1-t}\exp(M/k)}{2}.
}

Let
$
\phi(M,k,t):=1- \frac{k^{1-t}\exp(M/k)}{2}.
$
Note that for any fixed positive integers $M$ and $t$, $\phi(M,k,t)$ is increasing for $k\in (0,\infty)$.
Let $k_0:=\frac{2.4M}{t\log(M)}$  and 
\aln{eq1-1}
{
\psi(M,t):=2\cdot(2.4M)^{t-1}-(t\log(M))^{t-1}M^{\frac{t}{2.4}}.
}
Then $\psi(M,t)=2k_0^{t-1}\phi(M,k_0,t)(t\log(M))^{t-1}$. Thus, it suffices to show that $\psi(M,t)\ge 0$ when $t\ge 2$ and $M\ge t^3/2$.

When $M\ge t^3/2$, it is routine to check that $\psi(M,t)> 0$ holds for the case $2\le t\le 6$. 
In the following, we assume that $t\ge 7$.
Let 
$$\varphi_1(M):=2.4M-(2M)^{\frac13}\log(M)M^{\frac{1}{2.4}},$$
$$\varphi_2(M):=2^{\frac16}\cdot2.4M-(2M)^{\frac13}\log(M)M^{\frac{1}{2.4}+\frac{1}{2.4\cdot 6}}.$$
By applying fundamental techniques of Calculus, it can be easily verified that $\varphi_1(M),\varphi_2(M)> 0$ when $M> 0$. Then,
 \eqn{eq3-71}
 {
2\cdot(2.4M)^{t-1}
 &=&(2^{\frac16}\cdot2.4M)^6\cdot(2.4M)^{t-7} 
 \nonumber\\
 &>&
\left((2M)^{\frac13}\log(M)M^{\frac{1}{2.4}+\frac{1}{2.4\cdot 6}}\right)^6\cdot \left((2M)^{\frac13}\log(M)M^{\frac{1}{2.4}}\right)^{t-7}
\nonumber\\
&=&
(2M)^{\frac13(t-1)}(\log(M))^{t-1}M^{\frac{t}{2.4}}
\nonumber\\
&\ge&
(t\log(M))^{t-1}M^{\frac{t}{2.4}},
}
where the last inequality holds when $t\le (2M)^{\frac13}$.

Then the result follows from (\ref{eq1-1}) and (\ref{eq3-71}).
 \proofend

\noindent \textbf{Proof of Theorem~\ref{th4-1}.}
Let $M:=m-1$, $x:=\frac{M}{k}$ and $L$ be a $k$-assignment of $\hyh$ with $\alpha(\hyh,L)>0$.
Then by (\ref{eq6-1}),
\eqn{eq5-1}
{
\frac{P(\hyh, L)-P(\hyh,k)}
	{k^{n-3}\alpha(\hyh, L)}
	&>&  
		1
		-\sum_{i\ge 1}
	\frac{M^{2i-1}}{(2i-1)!}k^{-2i}
\nonumber \\
&=&1
-\frac{x}{M}
\left (\frac{\exp(x)-\exp(-x)}{2} \right )
\nonumber \\
&>&1
-\frac{x\exp(x)}{2M}.
}
Let $c:=0.844$ and
$$
 \phi(x,y):=1-\frac{x\exp(x)}{2y}.
$$
Then it can be verified that
\equ{eq5-54}
{
\phi(c\log (y),y)=
1 - \frac{c \log(y)}{2y^{1 - c}}
> \frac{0.002 \log(y)}{y^{1 - c}}
= \frac{0.002 \log(y)}{y^{0.156}}.
}
Note that for any fixed $y>0$, 
$\phi(x,y)$ is decreasing 
for $x\in (0,\infty)$.
Then (\ref{eq5-54}) implies that 
for any integer $M$ with $M\ge2$ and 
any real $x$ with 
$0<x\le c\log(M)$,
\equ{eq5-55}
{
\frac{P(\hyh, L)-P(\hyh,k)}
	{k^{n-3}\alpha(\hyh, L)}> \phi(x,M)
\ge \phi(c\log (M),M)>  \frac{0.002 \log(M)}{M^{0.156}}.
}
Since $x=\frac{M}{k}$ and $c=0.844$,
the result follows from  (\ref{eq5-1}) and (\ref{eq5-55}).
\proofend

\noindent\textbf{Proof of Theorem~\ref{th4-3}.}
Let $M:=m-1$, $x:=\frac{M}{k}$ and $L$ be a $k$-assignment of $\hyh$ with $\alpha(\hyh,L)>0$.
Then by (\ref{eq6-1}),
 \eqn{eq4-4}
{
 \frac{P(\hyh, L)-P(\hyh,k)}
 	{k^{n-r}\alpha(\hyh, L)}
 	&>&
 		1
 		-Mk^{-r+1}
 		+Mk^{-2r+4}
 		-\sum_{i\ge 1}
 	\frac{M^{2i-1}}{(2i-1)!}k^{-2i-2r+6}
\nonumber\\
& >&
 1-	\frac{x^{r-1}}{M^{r-2}}
-\frac{x^{2r-5}\exp(x)}{2M^{2r-5}}.
 }
Let
\equ{eq4-42}
{
\Psi(r):=
 1-	\frac{x^{r-1}}{M^{r-2}}-\frac{x^{2r-5}\exp(x)}{2M^{2r-5}}.
}
Note that as $\frac{x}{M}=\frac 1k< 1$, 
	$\Psi_1(r)<\Psi_1(r+1)$ 
	holds for all integers $r$ with $r\ge 4$.
	
	 Let $\phi(x,y):=2y^3-2yx^3-x^3\exp(x)$.
	Then $\Psi_1(4)=\frac{1}{2M^3}
	\phi(x,M)$.

	 Let $c:=\frac{1+(9/\exp(1))^{1/3}}{3}$. Then it can be verified that $0.830<c<0.831$. Thus for any $x\ge 0$, 
	  \eqn{eq4-6}
	  {
	  \exp(-x)\phi(x,\exp(cx))
	  &=&2\exp({(3c-1)x})-2x^3
	  \exp({(c-1)x})-x^3
	  \nonumber \\
	 & \ge&  
	  2\exp({(3c-1)x})-3x^3.
	  }
	  Let $\psi(x):=2\exp((3c-1)x)-3x^3$ and $x_0:=\frac{\log(M)}{c}$.
	 Then by (\ref{eq4-6}), 
	   \equ{eq4-7}
	  {
	  	\phi(x_0,M)
	  	=\phi(x_0,\exp(cx_0))
	  	\ge \exp(x_0)\psi(x_0) .
	  }
	 Observe that 
	  \equ{eq4-47}
	  {
	  	\left \{
	  \begin{array}{l}
	  \psi'(x)=2(3c-1)\exp((3c-1)x)-9x^2;\\
	  \psi''(x)=2(3c-1)^2\exp((3c-1)x)-18x;\\
	  \psi'''(x)=2(3c-1)^3\exp((3c-1)x)-18.
	  \end{array} 
	  \right.
	  }
	  
	  As $c=\frac{1+(9/\exp(1))^{1/3}}{3}$,
	  $x_1:=\frac{1}{3c-1}\log(\frac{9}{(3c-1)^3})$ is the unique zero of both $\psi'''(x)$ and $\psi''(x)$, and $\psi'''(x)\ge 0$ if and only if $x\ge x_1$.
	  Thus, $\psi''(x)\ge 0$ for all $x\ge 0$.
	  It further implies that $\psi'(x)\ge \psi'(0)=2(3c-1)>0$ for all $x\ge 0$ and thus
	  \aln{eq4-48}
	  {\psi(x_0)\ge \psi(0)+2(3c-1)x_0=2+2(3c-1)x_0.}
Note that for any fixed $y>0$, $\phi(x,y)$ 
 is decreasing for $x\in (0,\infty)$. Then by (\ref{eq4-7}) and (\ref{eq4-48}), for any integer $r$ with $r\ge 4$ and any real $x$ with $0<x\le x_0$,
 \eqn{eq4-8}
 {
 \Psi_1(r)\ge \Psi_1(4)&=&\frac{1}{2M^3}
 	\phi(x,M)
 	\nonumber\\
 	&\ge&  \frac{1}{2M^3}\phi(x_0,M)
 	\nonumber\\
 &\ge & \frac{1}{2M^3}\exp(x_0)(2+2(3c-1)x_0)
 \nonumber \\
 &=&M^{\frac 1c-3} \left (1+(3c-1)\frac{\log(M)}{c}\right )
 \nonumber \\
 &> & M^{-1.796}
 (1+1.796\log(M)),
 }
where the last inequality holds as $1/c>1.204$ and $M^{1/c}(3-1/c)>1.796M^{1.204}$ when $M\ge 2$.

Since $x=\frac{M}{k}$ and $c<0.831$, the result follows from (\ref{eq4-4}), (\ref{eq4-42}) and (\ref{eq4-8}).
 \proofend

\section*{Acknowledgement}

This research is supported by the Ministry of Education,
Singapore, under its Academic Research Tier 1 (RG19/22). Any opinions,
findings and conclusions or recommendations expressed in this
material are those of the authors and do not reflect the views of the
Ministry of Education, Singapore.
 
The first author would like to express her gratitude to National Institute of Education and Nanyang Technological University of Singapore for offering her Nanyang Technological University Research Scholarship during her PhD study.


\begin{thebibliography}{99}

	\bibitem{birk}
	G.D. Birkhoff, 
	A determinant formula for the number of 
	ways of coloring a map, \textit{Annal. Math.} \textbf{14}(1912), 42--46.
	

	
\bibitem{dong1} F.M. Dong and K.M. Koh, ``Foundations of the chromatic polynomial," in\textit{ the Handbook on the Tutte Polynomial and Related Topics}, Jo Ellis-Monaghan and Iain Moffatt (ed.), pp 232--266, CRC press, 2021.
	
	\bibitem{dong0} F.M. Dong, K.M. Koh and K.L. Teo, \textit{Chromatic Polynomials and Chromaticity of Graphs}, World Scientific, Singapore, 2005.
	
	\bibitem{Dong22b} F.M. Dong and M.Q. Zhang, Compare list-color functions 
	of uniform hypergraphs with their chromatic polynomials.
	\href{http://arxiv.org/abs/2212.02045}
	{At arXiv http://arxiv.org/abs/2212.02045.} 
	
	
\bibitem{Dong22} F.M. Dong and M.Q. Zhang, How large is $P(G,L)-P(G,k)$ for $k$-assignments $L$? 2022.
\href{http://arxiv.org/abs/2206.14536}{At arXiv http://arxiv.org/abs/2206.14536.}


\bibitem{erdos} P. Erd\H{o}s, A.L. Rubin and H. Taylor, Choosability in graphs, \textit{Congr. Numer.} \textbf{26} (1979), 125--157.


\bibitem{Jack15} B. Jackson, Chromatic polynomials, L.W. Beineke, R.J. Wilson (Eds.), \textit{Topics in Chromatic Graph Theory}, vol. 156, Cambridge University Press, 2015, pp. 56--72.

\bibitem{Kaul22} 
H. Kaul, A. Kumar, J.A. Mudrock, P. Rewers, P. Shin and K. To, On the List Color Function Threshold, arXiv:2202.03431 (preprint), 2022.

\bibitem{Kosto} A.V. Kostochka and A.F. Sidorenko, Problem Session of the Prachatice Con-ference on Graph Theory. In \textit{Fourth Czechoslovak Symposium on Combinatorics, Graphsand Complexity, Ann. Discrete Math,} vol. 51, p. 380, 1992.


\bibitem{rea1} R.C. Read, An introduction to chromatic polynomials, \textit{J. Combin. Theory} \textbf{4}(1968), 52--71.

\bibitem{rea2} R.C. Read and W.T. Tutte, Chromatic polynomials, in \textit{Selected Topics in Graph Theory 3}, Academic Press (1988), 15--42.

\bibitem{Thomassen} C. Thomassen, The chromatic polynomial and list colorings, \textit{J. Combin.Theory Ser. B} \textbf{99} (2009), 474--479.

\bibitem{Tomescu1998} 
I. Tomescu, Chromatic coefficients of linear uniform hypergraphs, 
{\it J. Combin. Theory Ser. B}
{\bf 260} (1998), 229–235.

\bibitem{Tomescu2009} 
I. Tomescu, Some properties of chromatic coefficients of linear uniform hypergraphs, 
{\it Graphs Combin.}
{\bf 25} (2009), 639–646.
	
	\bibitem{Trinks14} 
	M. Trinks, A note on a Broken-cycle Theorem for hypergraphs, \textit{Discuss. Math. Graph Theory} \textbf{34} (2014), 641--646.
	
	\bibitem{vizing} V.G. Vizing, Coloring the vertices of a graph in prescribed colors, \textit{Diskret. Analiz.no. 29, Metody Diskret. Anal. v Teorii Kodovi Skhem} \textbf{101} (1976), 3--10.
	
	\bibitem{wang20}
	W. Wang, J. Qian and Z. Yan, Colorings versus list colorings of uniform hypergraphs, \textit{J. Graph Theory} \textbf{95} (2020), 384--397.

\bibitem{wang17} W. Wang, J. Qian and Z. Yan, When does the list-coloring function of a graph equal its chromatic polynomial, \textit{J. Combin. Theory Ser. B} \textbf{122} (2017), 543--549.

\bibitem{Whitney1932} H. Whitney, A logical expansion in mathematics, \textit{Bull. Amer. Math. Soc.} \textbf{38} (1932) 572--579.

\bibitem{ruixue1} R.X. Zhang and F.M. Dong, Properties of chromatic polynomials of hypergraphs not held for chromatic polynomials of graphs, \textit{European Journal of Combinatorics}, \textbf{64}(2017), 138--151.

\bibitem{ruixue2} R.X. Zhang and F.M. Dong, Zero-free intervals of chromatic polynomials of hypergraphs, \textit{Discrete Mathematics} \textbf{343}(2020), 112--134.

\end{thebibliography}
\end{document}